\documentclass[10pt,a4paper]{amsart}
\usepackage{amsmath,amsthm,amssymb,amsfonts,mathrsfs}
\usepackage{geometry}
\newcommand\Sl{\textrm{SL}_2}

\newtheorem{Theorem}{Theorem}[section]
\newtheorem{Corollary}[Theorem]{Corollary}
\newtheorem{Lemma}[Theorem]{Lemma}
\newtheorem{Claim}[Theorem]{Claim}
\newtheorem*{NoNumberTheorem}{Theorem}
\newtheorem*{NoNumberProposition}{Proposition}

\theoremstyle{remark}
\newtheorem{Remark}{Remark}[section]

\theoremstyle{definition}
\newtheorem*{Acknowledgments}{Acknowledgments}

\newtheorem*{Pisot}{Pisot numbers}

\begin{document}
\title[Veech groups without parabolic elements]
{Veech groups without parabolic elements}

\author{Pascal Hubert, Erwan Lanneau}

\address{
Institut Math\'ematiques de Luminy (IML), UMR CNRS 6206 \\
163 avenue de Luminy, case 907\\
13288 Marseille cedex 9, France
}

\email{hubert@iml.univ-mrs.fr}

\address{
Max Planck Institute F\"ur Mathematics\\
Vivatsgasse 7, D-53111 Bonn, Germany
}

\email{lanneau@mpim-bonn.mpg.de}

\address{
Centre de Physique Th\'eorique (CPT), UMR CNRS 7061 \\
163 avenue de Luminy, case 907\\
13288 Marseille cedex 9, France
}

\email{lanneau@cpt.univ-mrs.fr}

\date{June 27, 2005}

\begin{abstract}
We prove that a translation surface which has two transverse
parabolic elements has totally real trace field. As a corollary,
non trivial Veech groups which have no parabolic elements do exist.

The proof follows Veech's viewpoint on Thurston's construction of
pseudo-Anosov diffeomorphisms.
\end{abstract}

\maketitle

\section{Introduction}

For a long time, it has been known that the ergodic properties of linear flows on a translation surface 
are strongly related to the behavior of its $\Sl (\mathbb R)$-orbit in the moduli space of 
holomorphic one forms
(see~\cite{MT}, \cite{ZO} for surveys of the literature on this subject).
The $\Sl( \mathbb R)$-orbit of a translation surface is called its Teichm\"uller disc. 
Its stabilizer under the action of  $\Sl( \mathbb R)$  is a Fuchsian group  called the Veech group.

In 1989, Veech proved that a translation surface whose stabilizer is a lattice has optimal 
dynamical properties: the directional flows are periodic or uniquely ergodic (see~\cite{Veech}). 
Since then, much effort has gone into the study of the geometry of Teichm\"uller discs 
(\cite{Ve3}, \cite{Vo}, \cite{Wa}).
Hubert and Schmidt~\cite{HS1,HS2} found the first examples of infinitely generated 
Veech groups. Just after that, McMullen~\cite{Mc1,Mc2} proved that, in genus 2, the existence 
of a pseudo-Anosov diffeomorphism in the affine group implies that the Veech group is a Fuchsian group 
of the first kind (which means that either it is a lattice or it is
infinitely generated; moreover McMullen proved that both cases
occur). See also~\cite{C} for related results.
McMullen's proof uses the existence of infinitely many parabolic
elements. By contrast, we will give examples with a very different
behavior in genus $g \geq 3$.

The trace field is a natural invariant of the Veech group. Thurston
proved that the trace of the derivative of any pseudo-Anosov diffeomorphisms is an algebraic
integer over $\mathbb Q$ with degree less that the dimension of the
Teichm\"uller space divided by $2$ (see~\cite{Thurston},~p.~427).
In~\cite{GJ} it was shown that a translation surface is a covering of
the torus ramified over one point if and only if its trace field
equals to $\mathbb Q$. Kenyon and Smillie~\cite{KS} gave 
a simple criterion ensuring this property: if the Veech group of a
translation surface contains a hyperbolic element whose trace belongs
to $\mathbb Q$, then this group is commensurable to $\Sl( \mathbb
Z)$. In fact they showed that the trace field is generated by the
trace of the derivative of any pseudo-Anosov diffeomorphisms. Moreover, if $K$ is the 
trace field of $(X,\omega)$ then the Veech group is commensurable to a subgroup 
of $\textrm{SL}_2(\mathcal O_K)$ where $\mathcal O_K$ is the ring of integers of $K$.

An interesting problem is to determine which
Fuchsian group can occur as affine group of some surface.
Up to now, there are no general methods to compute a Veech group.
To date there are two methods to produce pseudo-Anosov diffeomorphisms
in the coordinates of the flat surface. In the first 
one, due to Thurston,  a pseudo-Anosov diffeomorphism is obtained  as a product of two parabolic elements 
(see \cite{FLP}, \cite{Thurston}, \cite{Veech}). Veech computed the
first non trivial examples of affine groups by making calculations
with a pair of parabolic elements (see~\cite{Veech}).
Independently, a very  general 
construction of pseudo-Anosov diffeomorphisms 
was discovered by Veech~\cite{Ve1}. It is based on the Rauzy induction 
of interval exchange transformations (see also~\cite{Arnoux:Yoccoz} for  specific examples of such 
pseudo-Anosov diffeomorphisms, for any genus $g\geq 3$).
A simple consequence of our result is  that some pseudo-Anosov diffeomorphisms are not given 
by Thurston's construction (see \cite{Le} for another proof). In fact, we prove a stronger result:

\begin{Theorem}
\label{theo:1}
Let $(X,\omega)$ be a translation surface. Let us assume that the
Veech group $\textrm{SL}(X,\omega)$ contains two transverse
parabolic elements~\footnote{The surface $(X,\omega)$ is then called a {\it prelattice surface}~\cite{HS1} or 
a ``bouillabaisse surface'' in honor of John Hubbard's Lecture at the CIRM in July 2003.}.
Then the trace field
$$
\mathbb Q \left[\textrm{Trace}(A) \ , \ A \in
\textrm{SL}(X,\omega) \right]
$$
is totally real.
\end{Theorem}

When a pseudo-Anosov diffeomorphism $\phi$ acts linearly on a translation surface by the diagonal matrix 
$D\phi=\bigl[ \begin{smallmatrix} \lambda & 0 \\ 0 & \lambda^{-1}
\end{smallmatrix}\bigr]$ (with $\lambda^{-1} < 1 < \lambda$), its {\it
expansion factor} is $\lambda(\phi) = \lambda$ (and $1/\lambda$ is
the contraction factor). \medskip 

From Theorem~\ref{theo:1}, we draw the following results.

\begin{Theorem}
\label{theo:2}
Let $(X,\omega)$ be a translation surface endowed with a
pseudo-Anosov diffeomorphism~$\phi$ with expansion factor $\lambda$. Let
us assume that the field $\mathbb Q [\lambda + \lambda^{-1}]$ is not
totally real. Then $\textrm{SL}(X,\omega)$ does not contain any
parabolic elements.
\end{Theorem}

Arnoux and Yoccoz~\cite{Arnoux:Yoccoz} discovered a family
$\phi_n$, $n\geq 3$ of pseudo-Anosov diffeomorphisms with
expansion factor $\lambda_n=\lambda(\phi_n)$ the Pisot root of the
irreducible polynomial $P_n$ with
$$
P_n(X) = X^n-X^{n-1}-\dots-X-1.
$$
The pseudo-Anosov $\phi_n$ acts linearly on a genus $n$ surface (the corresponding
Abelian differential having two zeroes of order $n-1$). 

\begin{Corollary}
\label{cor:arnoux:yoccoz}
The Teichm\"uller disc stabilized by the Arnoux--Yoccoz
pseudo-Anosov $\phi_n$, $n\geq 3$  does not contain any
parabolic direction. \\
\noindent Therefore, for any genus $g\geq 3$, there exists a genus $g$ translation
surface such that its Veech group has (at least) one hyperbolic element and no parabolic
elements.
\end{Corollary}

\begin{Corollary} \label{Veech-surfaces}
The trace field of any Veech surface is totally real.
\end{Corollary}

\begin{Remark}
M\"oller proved Corollary \ref{Veech-surfaces} by very different methods (see \cite{Mo}).
\end{Remark}

\begin{Corollary}
\label{cor:pur:hyperbolic}
There exists a Veech group which is commensurable to a Fuchsian
group which only contains hyperbolic elements.
\end{Corollary}

\begin{Acknowledgments}
We thank John Hubbard for explaining to us Thurston's
construction. We also thank Chris Judge, Howard Masur, Curt McMullen, Thomas
Schmidt, Anton Zorich and the anonymous referees for helpful comments
on preliminary versions of this paper.

This work was done when the first author visited the second
at the Max-Planck-Institute f\"ur Mathematik in Bonn. We thank the
Institute for excellent working conditions.
\end{Acknowledgments}

\section{Background}

In order to establish notations and preparatory material, we
review basic notions concerning translation surfaces, affine
automorphisms groups and trace fields. We will end this section by
recalling  Veech's viewpoint on Thurston's construction. See
say~\cite{KS}, \cite{MT}, \cite{Mc1,Mc2}, \cite{Thurston}, \cite{Veech}
for more details; See also~\cite{Mc3,Mc4,Mc5}, for recent related
developments. For a general reference on Fuchsian groups, see~\cite{Katok}.

\subsection{Translation surfaces and affine diffeomorphisms groups}

A {\it translation surface} is a (real) genus $g$ surface with an atlas
such that all transition functions are translations. As usual, we
consider maximal atlases. These surfaces are precisely those given
by a Riemann surface $X$ and a holomorphic (non-null) one form
$\omega \in \Omega(X)$; see~\cite{MT} for a general reference
on translation surfaces and holomorphic one forms.

We denote by $X'$ the surface that arises from $X$ by deleting the
zeroes of the form $\omega$ on $X$. The translation structure
defines on $X'$ a Riemannian structure; we therefore have notions
of geodesic, length, angle, flow, measure... Orbits of the
directional flows meeting singularities are called separatrices.
Orbits of the flow going from a singularity to another one
(possibly the same) are called saddle connections.

Given any matrix $A\in \Sl(\mathbb R)$, we can post-compose the
coordinate functions of the charts of $(X,\omega)$ by $A$. 
One easily checks that this gives a new translation surface, denoted
by $A\cdot (X,\omega)$. We therefore get an $\Sl(\mathbb
R)$-action on these translation surfaces.

An {\it affine diffeomorphism $f : X \rightarrow X$} is a
homeomorphism of $X$ such that $f$ restricts to a diffeomorphism
on $X'$ of constant derivative. It is equivalent to say that $f$
restricts to an isomorphism of $X'$ which preserves the induced
affine structure given by $\omega$. Usually, one denotes by
$\textrm{Aff}(X,\omega)$ the group of orientation preserving
affine diffeomorphisms. The function which takes an affine
diffeomorphism $f$ to its derivative $D f$ gives a homomorphism
from $\textrm{Aff}(X,\omega)$ into $\Sl(\mathbb R)$. The image of
$\textrm{Aff}(X,\omega)$ is the {\it Veech group}
$\textrm{SL}(X,\omega)$ of the surface $(X,\omega)$ -- this is a
discrete subgroup and, when $X$ has genus greater than one, the
kernel of the homomorphism is finite.

One easily checks that the Veech group
$\textrm{SL}(X,\omega)$ is the $\Sl(\mathbb R)$-stabilizer of
$(X,\omega)$. Thus, for any matrix $A\in \Sl(\mathbb R)$, the
Veech group of $(X,\omega)$ and $A\cdot (X,\omega)=(Y,\alpha)$ are
conjugate in $\Sl(\mathbb R)$:
$$
\textrm{SL}(Y,\alpha) \ = \ A \cdot \textrm{SL}(X,\omega) \cdot
A^{-1}
$$

\subsection{Classification of affine diffeomorphisms}

There is a standard classification of the elements of $\Sl(\mathbb
R)$ into three types: elliptic, parabolic and hyperbolic. This
induces a classification of affine diffeomorphisms.

An affine diffeomorphism is respectively parabolic, elliptic or 
pseudo-Anosov if respectively $|\textrm{trace}(D f)| = 2$,
$|\textrm{trace}(D f)| < 2$ or $|\textrm{trace}(D f)| > 2$.

\begin{Remark}
\label{rk:finite:order}
If an elliptic element belongs to a Fuchsian group, its order is finite.
\end{Remark}

\begin{Remark}
\label{rk:parabolic:hyperbolic}
In a Fuchsian group, a parabolic direction (invariant direction of a 
parabolic element) is never fixed by a
hyperbolic element. More precisely, if a hyperbolic element $H$
fixes a parabolic direction of a parabolic element $P$ then one
can easily check that $H^nPH^{-n}$ converges to $\textrm{Id}$ as $n$ tends
to $+\infty$ (or $-\infty$), which is impossible in a discrete group.
\end{Remark}

\subsection{Cylinders decomposition and parabolic element}

A {\it cylinder} on $(X,\omega)$ is a maximal connected set of homotopic
simple closed geodesics. If the genus of $X$ is greater than one
then every cylinder is bounded by saddle connections. A cylinder
has a width (or circumference) $x$ and a height $y$. The {\it
modulus} of a cylinder is $\mu=y/x$. Veech~\cite{Veech} proved the
following:

\begin{NoNumberProposition}[Veech]
If a translation surface $(X,\omega)$ has a parabolic affine
diffeomorphism $f$, then there is a decomposition of $X$ into
metric cylinders parallel to the fixed direction of $D f$.
Furthermore, the moduli of the cylinders are commensurable (have
rational ratios).
\end{NoNumberProposition}

\begin{Remark}
\label{rk:puissance:twist}
In  the above proposition, up to take a power of the affine diffeomorphism, we 
can assume that $f$ acts as a power of  the affine  Dehn twist on each cylinder. 
Therefore the boundary  of each cylinder is fixed by $f$.
\end{Remark}
  
Conversely, a cylinder decomposition into cylinders of commensurable moduli
produces parabolic elements. Namely, the following holds:

\begin{NoNumberProposition}[Veech]
If $(X,\omega)$ has a decomposition into metric cylinders for the
horizontal direction, with commensurable moduli, then the Veech
group $\textrm{SL}(X,\omega)$ contains
$$
D f = \left( \begin{array}{cc}
1 & c \\
0 & 1
\end{array} \right)
$$
where $c$ is the least common multiple of the inverse of the moduli.
\end{NoNumberProposition}

\subsection{Trace fields}
\label{subsec-trace-field}

In this section we recall some general properties of the trace
field of a group; see~\cite{GJ}, \cite{KS}, \cite{Mc1,Mc2}.

The {\it trace field} of a group $\Gamma \subset \Sl(\mathbb R)$
is the subfield of $\mathbb R$ generated by tr$(A)$, $A\in
\Gamma$. One defines the trace field of a flat surface
$(X,\omega)$ to be the trace field of its Veech group
$\textrm{SL}(X,\omega) \subset \Sl(\mathbb R)$.

Let $(X,\omega)$ be a genus $g$ translation surface. Then the
following holds: \medskip

\noindent {\bf Theorem A.} (Kenyon,~Smillie). 
{\it The trace field of $(X,\omega)$ has degree at most $g$ over
$\mathbb Q$.

\noindent Assume that the affine diffeomorphisms group of $(X,\omega)$
contains a pseudo-Anosov element $f$ with expansion factor $\lambda$.
Then the trace field of $(X,\omega)$ is $\mathbb Q[\lambda +
\lambda^{-1}]$. } \medskip

One defines the {\it holonomy vectors} to be the integrals of
$\omega$ along the saddle connections. Let us denote
$\Lambda=\Lambda(\omega)$ the subgroup of $\mathbb R^2$ generated
by holonomy vectors
$$
\Lambda = \int_{H_1(X,\mathbb Z)} \omega
$$
Let $e_1,e_2 \in \Lambda$ be non-parallel vectors in $\mathbb
R^2$. One defines the {\it holonomy field} $k$ to be the smallest
subfield of $\mathbb R$ such that every element of
$\Lambda$ can be written $ae_1+be_2$ with $a,b \in k$. \medskip

\noindent {\bf Theorem B.} (Kenyon,~Smillie)
{\it The trace field of $(X,\omega)$ coincides with $k$. The space
$\Lambda \ \otimes \ \mathbb Q \subset \mathbb C$ is a
$2$-dimensional vector space over $k$.}

\medskip

See also~\cite{GJ} for a different approach of these notions. Note
that these results have been reproved in~\cite{Mc1,Mc2}.

\begin{Pisot}
An algebraic integer $\beta$ is a Pisot number if $\beta \in \mathbb
R$, $\beta > 1$ and all of its conjugates belong to the unit disc
$\mathbb D = \{ z\in \mathbb C,\ |z|<1 \}$.
\end{Pisot}

\subsection{Veech's viewpoint on Thurston's construction}
\label{sec:veech}

Let us recall the Thurston construction~\cite{Thurston}. We will follow the notations
of the paper of Veech~\cite{Veech}, section~$\S 9$.

Let $(Y,\alpha)$ be a translation surface with vertical and
horizontal parabolic directions. Up to  taking a power of the parabolic
elements, one can assume that the corresponding parabolic $P_{v}$
(resp $P_{h}$) is a multiple of the Dehn twist of each vertical
(resp horizontal) cylinder (see Remark~\ref{rk:puissance:twist}).

In these coordinates our two parabolic elements are
$$
P_h = \left( \begin{array}{cc}
1 & c \\
0 & 1
\end{array} \right) \qquad \textrm{and} \qquad
P_v = \left( \begin{array}{cc}
1 & 0 \\
d & 1
\end{array} \right)
$$

Without loss of the generality, we may assume that $c$ and $d$ are positive real numbers.

\begin{Claim} 
\label{claim:tracefield}
Let $t=cd>0$; then the trace
field of $\textrm{SL}(Y,\alpha)$ is $\mathbb Q[t]$.
\end{Claim}

\begin{proof}[Proof of Claim~\ref{claim:tracefield}]
The matrix  $P_h P_v$ has trace $2+t > 2$, thus
this is a hyperbolic element and, following~\cite{KS} (see section~\ref{subsec-trace-field}, 
Theorem~A), the trace field of $\textrm{SL}(Y,\alpha)$ is $\mathbb Q[t]$. So the claim is proven.
\end{proof}

Let us denote by $H_i$, $1 \leq i \leq r$ and $V_j$, $1 \leq j
\leq s$ the horizontal and vertical cylinders. Let us denote the
width and heights of $H_i$ and $V_j$ respectively by $(x_i,y_i)$
and $(\eta_j,\xi_j)$. We insist that the first coordinate is the width and the second 
one is the height even for vertical cylinders.

With these notations, let $E$ be the $r
\times s$ integer matrix whose entry $E_{i,j}$ is the number of
rectangles $(\xi_j\times y_i)$ in the intersection $H_i \cap V_j$.
All of these rectangles have width $y_j$ and heights $\xi_i$. In
other words, $E_{i,j}$ is the intersection number of the core
curves of the cylinders $H_i$ and $V_j$.

Let us introduce the following notations of linear algebra:
$\overrightarrow{x}=(x_1,\dots,x_r)$,
$\overrightarrow{y}=(y_1,\dots,y_r)$,
$\overrightarrow{\xi}=(\xi_1,\dots,\xi_s)$ and
$\overrightarrow{\eta}=(\eta_1,\dots,\eta_s)$. Then one can
summarize the above discussion by the matrix relations:
\begin{equation}
\label{eq:1}
\left\{ \begin{array}{c}
\overrightarrow{x} = E \overrightarrow{\xi} \\
\overrightarrow{\eta} = {}^{\textrm t}E \overrightarrow{y}
\end{array} \right.
\end{equation}
The moduli of the vertical cylinder $V_j$ (resp horizontal
cylinder $H_i$) is commensurable with $d$ (resp with $c$). More
precisely, there exist integers $m_i$, $1 \leq i \leq r$, and
$n_j$, $1 \leq j \leq s$, such that
\begin{equation}
\label{eq:veech}
\left\{ \begin{array}{c}
m_i x_i = c y_i \\
n_j \eta_j = d \xi_j
\end{array} \right.
\end{equation}
Let us denote by $D_m=\textrm{Diag}(m_1,\dots,m_r)$ and
$D_n=\textrm{Diag}(n_1,\dots,n_s)$ the diagonal matrices. Then the
above equation~(\ref{eq:veech}) becomes:
\begin{equation}
\label{eq:2}
\left\{ \begin{array}{c}
D_m \overrightarrow{x} = c \overrightarrow{y} \\
D_n \overrightarrow{\eta} = d \overrightarrow{\xi}
\end{array} \right.
\end{equation}
From equations~(\ref{eq:1}) and~(\ref{eq:2}) one gets the
following new one:
\begin{equation*}
\left\{ \begin{array}{c}
E D_n \overrightarrow{\eta} = d \overrightarrow{x} \\
{}^{t}E D_m \overrightarrow{x} = c \overrightarrow{\eta}
\end{array} \right.
\end{equation*}
and therefore we deduce:
\begin{equation}
\label{eq:4}
\left\{ \begin{array}{c}

E D_n {}^{t}E D_m \overrightarrow{x} = cd \overrightarrow{x} \\
{}^{t}E D_m E D_n \overrightarrow{\eta} = cd \overrightarrow{\eta}
\end{array} \right.
\end{equation}
Now, in order to follow Veech's notations, let us introduce
the two matrices $F_n=E D_n$ and $F_m= {}^{t}E D_m$. As remarked
in~\cite{Veech}, the matrices $F_n F_m$ and $F_m F_n$ have a power
with positive entries (see \cite{HuLe} Appendix C for a proof).
The vector $\overrightarrow{x}$ is a non negative eigenvector of the 
Perron--Frobenius matrix $F_n F_m$, therefore $t=cd>0$ is the
unique Perron--Frobenius eigenvalue of $F_n F_m$ (the same is true for
the Perron--Frobenius matrix $F_m F_n$).
Thus, up to renormalization of the area of the surface, 
the coordinates of
$\overrightarrow{x}$ and $\overrightarrow{\eta}$ belong to
$\mathbb Q[t]$ (see section~\ref{subsec-trace-field}, Theorem~B).

\vskip 5mm

Now we have all necessary tools to prove the announced results.

\section{Proofs}

We first prove Theorem~\ref{theo:2} assuming Theorem~\ref{theo:1}.

\begin{proof}[Proof of Theorem~\ref{theo:2}]
Let us assume that there is a parabolic element $P$ in
$\textrm{SL}(X,\omega)$. Let us denote by $H$ the derivative of
the pseudo-Anosov $\phi$. Then the conjugate $H P H^{-1}$ is
another parabolic element in $\textrm{SL}(X,\omega)$. Let $x\in
\partial \mathbb H$ be the fixed point of $P$. Thus, $H(x)$ is a
fixed point of $H P H^{-1}$. But by
Remark~\ref{rk:parabolic:hyperbolic}, $H(x) \not = x$, then
$H P H^{-1}\in \textrm{SL}(X,\omega)$ is certainly a parabolic
element transverse to the parabolic $P$. Therefore
Theorem~\ref{theo:1} applies.
\end{proof}

\begin{proof}[Proof of Theorem~\ref{theo:1}]
Let us assume that the surface $(X,\omega)$ has two transverse
parabolic elements. By a standard argument, one can find a matrix
$A\in\textrm{SL}_2(\mathbb R)$ which sends the two invariant
directions of our parabolic elements into horizontal and vertical
direction. The Veech group
$$
\textrm{SL}(Y,\alpha)=A \cdot \textrm{SL}(X,\omega) \cdot A^{-1}
$$
possesses the same trace field as $\textrm{SL}(X,\omega)$.

Now up to taking a power of the parabolic element, one can assume
that they are a multiple of the Dehn twist on each vertical (resp
horizontal) cylinder (see Remark~\ref{rk:puissance:twist}). Thus
one can apply Veech's viewpoint on Thurston's construction,
section~\ref{sec:veech}. In particular we follow the notations
introduced in that section.

Recall that the trace field of  $\textrm{SL}(Y,\alpha)$ is
$\mathbb Q[t]$ (see Claim \ref{claim:tracefield}).
 Now let us prove that $\mathbb Q[t]$
is totally real.

\vskip 5mm

Let $\sigma$ be an embedding of $\mathbb Q[t]$ into $ \mathbb C$  and 
$t'=\sigma(t) \in \mathbb C$ be a conjugate of
$t$. Applying $\sigma$ to the first part of
equation~(\ref{eq:4}): $F_n F_m \overrightarrow{x}=t
\overrightarrow{x}$ and recalling that $F_n F_m$ is an {\it
integer matrix}, one gets
\begin{equation}
\label{eq:new:conjugate}
F_n F_m \sigma(\overrightarrow{x})=t' \sigma(\overrightarrow{x})
\end{equation}
Now, let us denote by
$D_{\sqrt{m}}=\textrm{Diag}(\sqrt{m_1},\dots,\sqrt{m_r})$ and
$D_{\sqrt{n}}=\textrm{Diag}(\sqrt{n_1},\dots,\sqrt{n_s})$ the
diagonal matrices. Then
$$
F_n F_m = E D_n {}^{t}E D_m = E D_{\sqrt{n}} \ D_{\sqrt{n}}
{}^{t}E \ D_{\sqrt{m}} D_{\sqrt{m}} \ = \ E D_{\sqrt{n}} \
{}^{t}(E D_{\sqrt{n}} ) \ D_{\sqrt{m}} D_{\sqrt{m}}
$$
Let us set $A=E D_{\sqrt{n}}$. Substituting this into the last
equation, yields:
$$
F_n F_m = A {}^{t} A D_{\sqrt{m}} D_{\sqrt{m}}
$$
Letting $M =  D_{\sqrt{m}} \ A$, it becomes:
\begin{multline}
\label{eq:semblable} F_n F_m \ = \ D_{\sqrt{m}}^{-1} D_{\sqrt{m}}
\ A {}^{t} A D_{\sqrt{m}} D_{\sqrt{m}}
 \ = \
\\
 = \ D_{\sqrt{m}}^{-1} D_{\sqrt{m}} \ A \ {}^{t}(D_{\sqrt{m}} \ A) \ D_{\sqrt{m}}
\ = \qquad \qquad
\\
= \ D_{\sqrt{m}}^{-1} M {}^{t}M \ D_{\sqrt{m}} \qquad \qquad
\qquad \qquad \qquad
\end{multline}
Now equation~(\ref{eq:new:conjugate}) and the fact that
$\sigma(\overrightarrow{x}) \not = \overrightarrow{0}$ imply
that $t'$ is an eigenvalue of $F_n F_m$. But by
equation~(\ref{eq:semblable}), the two matrices $F_n F_m$ and $M
{}^{t}M$ are similar, they thus have the same eigenvalues. But $M
{}^{t}M$ is symmetric, thus all of its eigenvalues are real, and
so $t' \in \mathbb R$.

Finally the trace field $\mathbb Q[t]$ of $(Y,\alpha)$, and
that of $(X,\omega)$, is totally real. Theorem~\ref{theo:1} is
proved.
\end{proof}

\begin{proof}[Proof of Corollary~\ref{cor:arnoux:yoccoz}]
Let $n\geq 3$ be any odd integer. We denote by $(X_n,\omega_n)$ a
flat surface in the Teichm\"uller disc stabilized by the
Arnoux--Yoccoz pseudo-Anosov $\phi_n$. By Theorem~A 
(see section~\ref{subsec-trace-field}), the trace field of
$(X_n,\omega_n)$ is $\mathbb Q[\lambda_n + \lambda_n^{-1}]$.

\begin{Claim}
\label{claim:racine}
The polynomial $X^n - X^{n-1} - \dots - 1$ has $2$ real roots if $n$
is even and $1$ if $n$ is odd.
\end{Claim}

\begin{proof}[Proof of the Claim]
Following~\cite{Arnoux:Yoccoz} we introduce the polynomial
$$
Q_n(X) = (X^n - X^{n-1} - \dots - 1)(X-1) = X^{n+1}-2X^n + 1
$$
One can directly check, by calculating $Q'_n$, that $Q_n(X)$ has $2$
real roots if $n$ is odd and $3$ if $n$ is even. This proves the claim.
\end{proof}

Above Claim~\ref{claim:racine} asserts that $\mathbb Q[\lambda_n]$ is
not totally real. Recall that $\lambda_n$ is a Pisot number
(see~\cite{Arnoux:Yoccoz}). Applying next Lemma~\ref{lm:pisotnumber}
we get that $\mathbb Q[\lambda_n + \lambda_n^{-1}]$ is
not totally real.

Thus Corollary~\ref{cor:arnoux:yoccoz} follows from Theorem~\ref{theo:2}.
\end{proof}

\begin{Lemma}
\label{lm:pisotnumber}
Let $\beta$ be any Pisot number. Let us assume that $\mathbb Q[\beta]$
is not totally real. Then the field $\mathbb Q[\beta + \beta^{-1}]$ is not
totally real. 
\end{Lemma}
 
\begin{proof}[Proof of Lemma~\ref{lm:pisotnumber}]
Let $\delta$ be a conjugate of $\beta$ which is not real. Galois
theory ensures that there is a field homomorphism $\chi : \mathbb
Q[\beta] \rightarrow \mathbb Q[\delta]$. The complex number $\chi(\beta
+ \beta^{-1})=\delta + \delta^{-1}$ is a conjugate of $\beta + \beta^{-1}$.
It is enough to show that $\delta + \delta^{-1}$ is not real to 
prove that $\mathbb Q[\beta +  \beta^{-1}]$ is not totally real. Writing 
$\delta = \rho \displaystyle e^{i\theta}$ (with $\sin(\theta) \not = 0$), 
we have $\Im m(\delta + \delta^{-1}) = (\rho - \rho^{-1})\sin( \theta)$.
As $\beta$ is a Pisot number, $\rho =\vert \delta \vert < 1$. Therefore 
$\delta + \delta^{-1}$ is not real.
So Lemma~\ref{lm:pisotnumber} is proven.
\end{proof}

\begin{proof}[Proof of Corollary \ref{Veech-surfaces}]
On a Veech surface, the direction of every saddle connection is a parabolic direction. 
There are thus at least two transverse parabolic elements in the Veech group and 
Theorem~\ref{theo:1} applies.
\end{proof}

\begin{proof}[Proof of Corollary~\ref{cor:pur:hyperbolic}]
Let $(X,\omega)$ be any genus $g\geq 3$ translation surface whose
Veech group only contains hyperbolic and elliptic elements. Any
elliptic element in $\textrm{SL}(X,\omega)$ is conjugate in
$\Sl(\mathbb R)$ to a rotation. As a rotation preserves the
underlying complex structure of the Riemann surface $X$, it is an
automorphism of a genus $g$ Riemann surface. Therefore, by a
Theorem of Hurwitz, the order of any elliptic element belonging to
$\textrm{SL}(X,\omega)$ is bounded by $84(g-1)$, see
say~\cite{Farkas:Kra} $\S 5$ p.$242$.

Now we recall a Theorem of Purzitsky on Fuchsian groups 
(see~\cite{Pu} Theorem~$7$ p.$241$): 

\begin{NoNumberTheorem}[Purzitsky]
Let $\Gamma$ be a Fuchsian group. Then $\Gamma$ contains a finite index
subgroup without elliptic elements if and only if
there exists a constant $N$ such that the order of any elliptic
element of $\Gamma$ is less than $N$.
\end{NoNumberTheorem}

Now recalling that any elliptic element belonging to a Fuchsian group
has finite order (see Remark~\ref{rk:finite:order}),
Corollary~\ref{cor:pur:hyperbolic} follows from Purzitsky's 
Theorem taking $\Gamma=\textrm{SL}(X,\omega)$.
\end{proof}


\end{document}